%% file: braidstar123004.tex
\documentclass[12pt]{article}

\usepackage{amsmath}
\usepackage{amsthm}
\usepackage[dvips]{graphics}
\usepackage{epsfig}
\usepackage{color}

\newtheorem{thm}{Theorem}
\newtheorem{lemma}[thm]{Lemma}

\newtheorem{cor}[thm]{Corollary}
\newtheorem*{ex}{Example}

\theoremstyle{definition}

\title{Stellar Braiding}
\author{Margaret Doig \\ University of Notre Dame \\ mdoig@nd.edu}

\begin{document}

\maketitle

\begin{abstract}
In this paper, we exhibit an explicit one-dimensional deformation
retract $D_k(S_n)$ of the unordered $k$-point configuration space
$U_k^{top}(S_n)$ for any star $S_n$.  These spaces have recently
been studied by Abrams and Ghrist.  We use this retract to compute
an explicit set of free generators $\beta_k$ of the corresponding
braid group $B_k(S_n,c_k)$.  In particular, we show that the
natural map $i_{k*}:B_{k-1}(S_n,c_{k-1}) \hookrightarrow
B_k(S_n,c_k)$ sends $\beta_{k-1}$ to $\beta_k$ injectively.
\end{abstract}

\section{Introduction}

We investigate the braid group of a star, or the fundamental group
of its configuration space.  The \emph{unordered k-point
configuration space} of a space $X$ is the set of $k$-element
subsets, or \[U^{top}_k(X)=\{c \subset X : |c|=k\}.\]  This space
is topologized as the quotient space of the map \[p: X^k -
\Delta(X) \rightarrow U^{top}_k(X)\] such that $p(x_1, x_2,
\cdots, x_k)=\{x_1, x_2, \cdots, x_k\}$ where the \emph{fat
diagonal} of $X$ is \[\Delta(X)=\{(x_1, x_2, \cdots, x_k) \in X^k
: \exists i \neq j~such~that~x_i = x_j\}.\] The product topology
on $X^k-\Delta(X)$ induces a topology on $U^{top}_k(S_n)$.

The \emph{k-point braid group} of a star, $B_k(S_n,c_k)$, is the
fundamental group of the configuration space
$\pi_1(U_k^{top}(S_n),c_k)$ for some arbitrary basepoint $c_k$.  A
\emph{star} $S_n$ is a graph on $n+1$ vertices with one vertex
$v_0$ of degree $n$ and $n$ vertices $v_1, v_2, \cdots, v_n$ of
degree one.  On $S_n$, we use the metric
\[d(x,y)=\kappa \rho (x,y)~\forall x,y \in S_n\] where $\rho$ is
the standard simplicial metric on $S_n$ and $\kappa$ is a fixed
constant greater than or equal to $k-1$. Therefore
$d(v_0,v_i)=\kappa$ if $i > 0$ and $d(v_i,v_j)=2 \kappa$ if $1
\leq i < j \leq n$.  This induces a Hausdorff metric on the
collection of non-empty closed subsets of $S_n$:
\[d(A,B)=\sup_{\substack{a \in A \\ b \in B}} (d(a,B),d(A,b)).\]

This problem arises in an overlap of topology and graph theory.
Graphs treated as topological spaces are well understood, as are
configuration spaces of Euclidean spaces.  Recently, Abrams and
Ghrist introduced the notion of the $k$-point braid group of a
graph (see~\cite{factory} for an overview of the area).  This
group is much less clearly understood but is in many ways
analogous to Artin's classical braid group, the braid group of the
real plane. In the case of a star $S$, Ghrist has shown that the
$k$-point braid group $B_k(S,c_0)$ is a free group. He has also
shown that the associated $k$-point configuration space
$U_k^{top}(\Gamma)$ for any graph $\Gamma$ is an aspherical space
(see~\cite{ghrobotics} Theorem 3.1). Further, he has shown that
the natural map $B_{k-1}(\Gamma,c_0) \hookrightarrow
B_k(\Gamma,c_0')$ is injective if $\Gamma$ is a tree (Lemma 2.2).

In Section~\ref{polyhedron}, we exhibit a one-dimensional
polyhedron $D_k(S_n)$.  We show in Section~\ref{retract} that the
polyhedron is a deformation retract of $U_k^{top}(S_n)$ and name
it the \emph{k-spine} of $S_n$. We then use $D_k(S_n)$ to generate
$B_k(S_n,c_k)$ directly.  Finally, we show that the natural
simplicial map $i_k:D_{k-1}(S_n) \hookrightarrow D_k(S_n)$ induces
a monomorphism $i_{k*}:B_{k-1}(S_n,c_{k-1}) \hookrightarrow
B_k(S_n,c_k)$ mapping a basis $\beta_{k-1}$ of
$B_{k-1}(S_n,c_{k-1})$ into a basis $\beta_k$.

\begin{figure}
\begin{center}
\input{vertices.pstex_t} \caption{Vertices in $D_k(S_n)$ are (a)
type I if $d(c,v_0)=0$ and (b) type II if $d(c,v_0)=1/2$
(\emph{N.B.:} in (b), $A_1(c)=A_5(c)=A_6(c)=\emptyset$ and so are
not arms).} \label{vertices}
\end{center}
\end{figure}
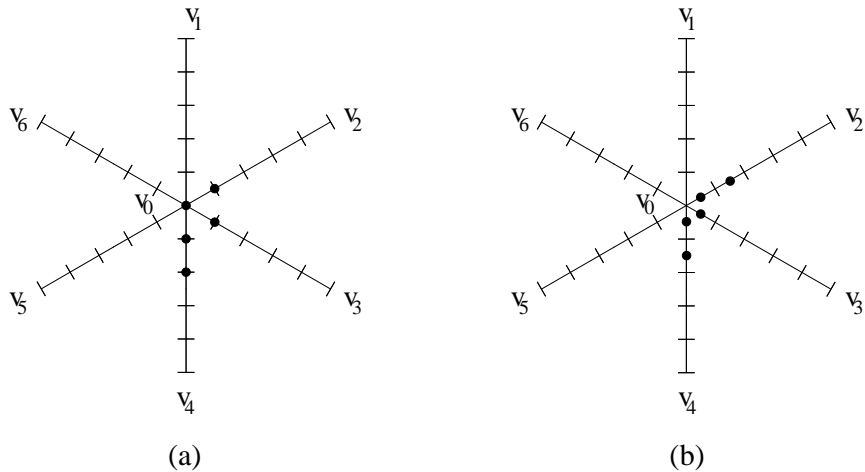

\section{The space $D_k(S_n)$}\label{polyhedron}
We construct a one-dimensional polyhedron $D_k(S_n)$ within
$U_k^{top}(S_n)$.  Later we will show that it is a deformation
retract of $U_k^{top}(S_n)$. Number the vertices $\{v_0, v_1,
\cdots, v_n\}$ of $S_n$ so that $deg(v_0)=n$.  The geodesic
segment in $S_n$ joining any two points $x,~y$ is written $[x,y]$.
We say a finite subset $c$ of $S_n$ is a \emph{chain} if, for each
$x,y \in c$, we have: \[|c \cap [x,y]|=1+d(x,y).\] Roughly
speaking, each point of $c$ is of distance one from its neighbors.
We also define: \[A_i(c)=c \cap [v_0, v_i]\] for $c \subset S_n$
and call all nonempty $A_i(c)$ the \emph{arms} of $c$.  Then we
let a configuration $c$ be \emph{regular} if there is an arm
$A_m(c)$ called a \emph{governing arm} where $A_m(c) \cup A_i(c)$
is a chain for all $i$ and there are at least two arms of $c$.
Set
\[D_k(S_n)=\{c:c = regular\}.\] That is, for any $c \in D_k(S_n)$,
there is a governing arm $A_m(c)$ whose union with every other arm
is a chain.

If the governing arm is unique, then $0< d(c,v_0) = d(A_m(c),v_0)
< 1/2$; moreover, each arm is a chain, and $d(A_i(c),v_0)=1-
d(c,v_0)$ for every other arm $A_i(c)$.  On the other hand,
$A_m(c)$ is not unique exactly when $d(A_i(c),v_0)=0$ or 1/2 for
each arm $A_i(c)$.  In either case, each $c \in D_k(S_n)$ is
uniquely determined by the numbers $|A_i(c)|$, by an integer $m$
where $A_m(c)$ is a governing arm, and by the value $d(c,v_0)$. We
can now prove the following about $D_k(S_n)$:

\begin{lemma}\label{graph}
The space $D_k(S_n)$ is a one-dimensional polyhedron, i.e., a graph.
\end{lemma}
\begin{proof}
We define a set of vertices of $D_k(S_n)$ by: \[V(D_k(S_n)) =\{c
\in D_k(S_n): d(c,v_0)=0~or~1/2\}.\]

There are two types of vertices in $V(D_k(S_n))$.  A type I vertex
is an element $c \in V(D_k(S_n))$ such that $v_0 \in c$, as in
Figure~\ref{vertices}a. There is a bicorrespondence between
vertices of this type and distributions of $k-1$ identical
particles into the $n$ disjoint edges $(v_0,v_i]$, and so there
are $\binom{n+k-2}{n-1}$ of these type I vertices.  On the other
hand, a vertex $c$ is of type II if $d(c,v_0)=1/2$ like
Figure~\ref{moves}b. There is one of these for each way of
distributing $k$ particles in the $n$ edges $(v_0,v_i]$ without
putting all the particles on the same edge. Thus, there are
$\binom{n+k-1}{n-1}-n$ vertices of this type.

We next define the set of edges of $D_k(S_n)$.  Take two vertices
$c_0$ and $c_1$ such that $d(c_0,c_1)=1/2$.  They will determine
an edge of $D_k(S_n)$ denoted $[c_0,c_1]$.  Note that exactly one
of $c_0$ and $c_1$ has type I, say $c_0$.  Then $c_1$ has type I.
There is a unique $m$ such that $|A_m(c_0)|=|A_m(c_1)|$.  It
follows that $|A_i(c_0)|=|A_i(c_1)|+1$ for $i \neq m$.
Additionally, $d(A_i(c_0),A_i(c_1))=1/2$ for each arm $A_i(c_0)$.
We may now define the edge $[c_0,c_1]$ of $D_k(S_n)$ by:
\[[c_0,c_1]=\{c \in D_k(S_n):d(c,c_0) + d(c,c_1) = 1/2\}.\]  Note
that a configuration $c \in D_k(S_n)$ is in $[c_0,c_1]$ exactly
when $d(A_i(c),A_i(c_0)) \leq 1/2$ and $d(A_i(c),A_i(c_1))\leq
1/2$.  Note that $A_m(c)$ is the unique governing arm of $c$ when
$c \neq c_0$ or $c_1$ and that $|A_i(c)|=|A_i(c_1)|$ for $i \neq
m$. The set of edges of $D_k(S_n)$ is: \[E(D_k(S_n))=\{[c_0,c_1] :
c_0,c_1 \in V(D_k(S_n))~and~d(c_0,c_1)=1/2\}.\]

To complete the proof, we claim that each configuration in
$D_k(S_n) / V(D_k(S_n))$ lies in exactly one edge $[c_0,c_1]$ and
that each edge $[c_0,c_1]$ is homeomorphic to $[0,1]$.  To see
this, let $c \in D_k(S_n)/ V(D_k(S_n))$.  Recall each
configuration is uniquely determined by the numbers $|A_i(c)|$, by
an integer $m$ where $A_m(c)$ is a governing arm, and by the value
$d(c,v_0)$.  Then let  $c_0$ be the type I vertex uniquely
determined by $|A_m(c_0)|$=$|A_m(c)|$, by $|A_i(c_0)|-1=|A_i(c)|$
for all $i \neq m$.  let $c_1$ be the type II vertex uniquely
determined by $|A_i(c_1)|=|A_i(c)|$ for all $i$.  Set
$p=d(c,v_0)$.  Then $d(c,c_0)=p \leq 1/2$ and $d(c,c_1)=1/2-p$.
Therefore $c \in [c_0,c_1]$. Additionally, $c$ is uniquely
determined by $c_0$,$c_1$, and $p$.  A homeomorphism \[h:
[c_0,c_1] \rightarrow [0,1]\] is then given by the rule
$h(c)=2d(c,v_0)$.
\end{proof}

\section{A deformation retraction}\label{retract}

\begin{figure}
\begin{center}
\input{moves.pstex_t}
\caption{Deformation retract $R$ taking $U_5^{top}(S_6)$ to
$D_5(S_6)$, mapping (a) configuration $c$ to a type I vertex $c'$
with $d(c',v_0)=0$ and (b) configuration $c$ to a non-vertex point
$c'$ with governing arm $A_4(c')$, that is, with
$d(A_4(c'),v_0)=p<1/2$ and $d(A_i(c'),v_0)=1-p$ for all $i$ where
$A_i(c) \neq \emptyset$ and $A_i(c')$ is an arm.} \label{moves}
\end{center}
\end{figure}
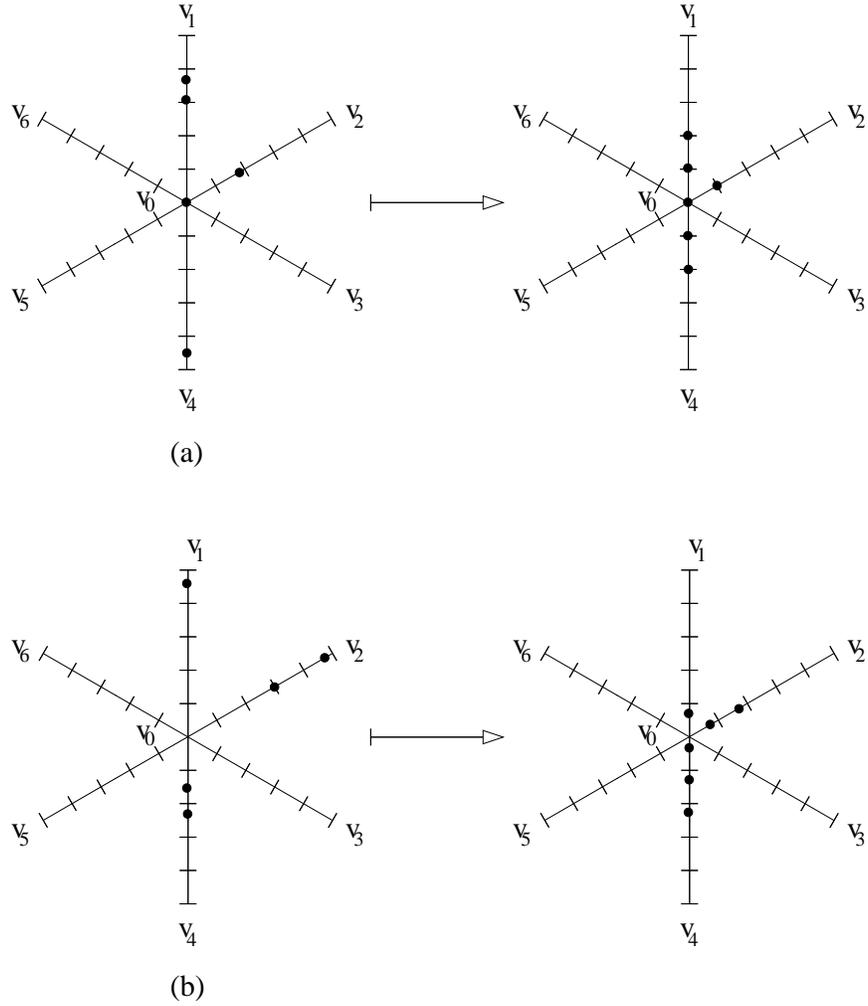

The construction of $D_k(S_n)$ suggests a deformation retraction:

\begin{thm}\label{def}
$D_k(S_n)$ is a deformation retract of $U_k^{top}(S_n)$.
\end{thm}
\begin{proof}
Let $c \in U_k^{top}(S_n)$.  Without loss of generality, assume
$d(A_1(c),v_0) \leq d(A_2(c),v_0) \leq d(A_i(c),v_0)$ for all $i
\geq 3$.  Set $p=d(A_1(c),v_0)$ and $q=d(A_2(c),v_0)$.  We are
going to define a homotopy $R: U_k^{top}(S_n) \times [0,1]
\rightarrow U_k^{top}(S_n)$ such that $R(c,0)=c$ and $R(c,1)=c'
\in D_k(S_n)$. If $|A_i(c)|=0$ for all $i$, a regular
configuration $c'$ is uniquely defined by $|A_i(c')|=|A_i(c)|$ and
by $d(A_i(c'),v_0)=0$ for all $i$. On the other hand, if
$|A_i(c)|\neq 0$ for some $i$, then $|A_i(c)|\neq 0$ for all $i$,
and then a regular configuration $c'$ is uniquely defined by
$|A_i(c')|=|A_i(c)|$ and by $d(A_1(c'),v_0)=\frac{p}{p+q}$ and
$d(A_i(c'),v_0) =\frac{q}{p+q}$ for all $i > 1$.  There is a
natural ordering on the sets $A_i(c)$ and $A_i(c')$ given by
$d(x_{ij},v_0)$ for $x_{ij} \in A_i(c)$ or $A_i(c')$.  Then
$A_i(c)$ and $A_i(c')$ correspond to points in $[0,k-1] \times
[0,k-1]\times \cdots \times [0,k-1]$ below the fat diagonal, and
so the straight line segment between them does not cross the fat
diagonal and gives a homotopy $h_i$.  These coalesce to provide a
homotopy $R: U_k^{top}(S_n) \times [0,1] \rightarrow
U^{top}_k(S_n)$. See Figure~\ref{moves}.
\end{proof}

We may now calculate $B_k(S_n,c_0)$ via $\pi_1(D_k(S_n),c_0)$. Let
$c_0$ be the unique vertex contained in $[v_0,v_1]$.  Ghrist gives
the order of $\pi_1(C_k(S_n),c_0)$, the fundamental group of the
ordered configuration space, i.e., of $U_k^{top}(S_n)$ where each
element is an ordered $n$-tuple rather than a
set~\cite{ghrobotics}. He uses an inductive proof; we demonstrate
the order follows as a corollary of Theorem~\ref{def}.

\begin{cor}\label{group}
The braid group $B_k(S_n,c_0)$ is free on $\varphi(k,n)$
generators where \[\varphi(k,n)=1 + (n-1)
\binom{n+k-1}{n-1}-\binom{n+k-1}{n-1}.\]
\end{cor}
\begin{proof}
Note that $D_k(S_n)$ is one-dimensional and so has a free group as
fundamental group.  The rank can be computed from the Euler
characteristic $\chi$ of $D_k(S_n)$.  We can count the edges in
$E(D_k(S_n))$.  Since each edge contains exactly one vertex of
type I, we need merely add the degrees of the type I vertices.
Each of the $n$ vertices $c_0$ which are contained entirely in one
arm $A_m(c_0)$ have degree $n-1$ because each is adjacent to the
$n-1$ vertices $c_1$ where $|A_m(c_1)| = k-1$ and $|A_i(c_1)|=1$
for $i \neq m$. On the other hand, the remaining
$\binom{n+k-1}{n-1}-n$ type I vertices $c_0$ have degree $n$
because each is adjacent to $n$ different $c_1$ where
$|A_m(c_1)|=|A_m(c_0)|$ for some $m$ and $|A_i(c_1)|=|A_i(c_0)|-1$
for all $i \neq m$.  Thus
\[|E(D_k(S_n))|=n\binom{n+k-1}{n-1}-n.\]  We can also count the
vertices since
\[|V(D_k(S_n))|=\binom{n+k-2}{n-1}+\binom{n+k-1}{n-1}-n.\] Then
\[rank(B_k(S_n))=1-\chi=1+(n-1)
\binom{n+k-2}{n-1}-\binom{n+k-1}{n-1}.\]
\end{proof}

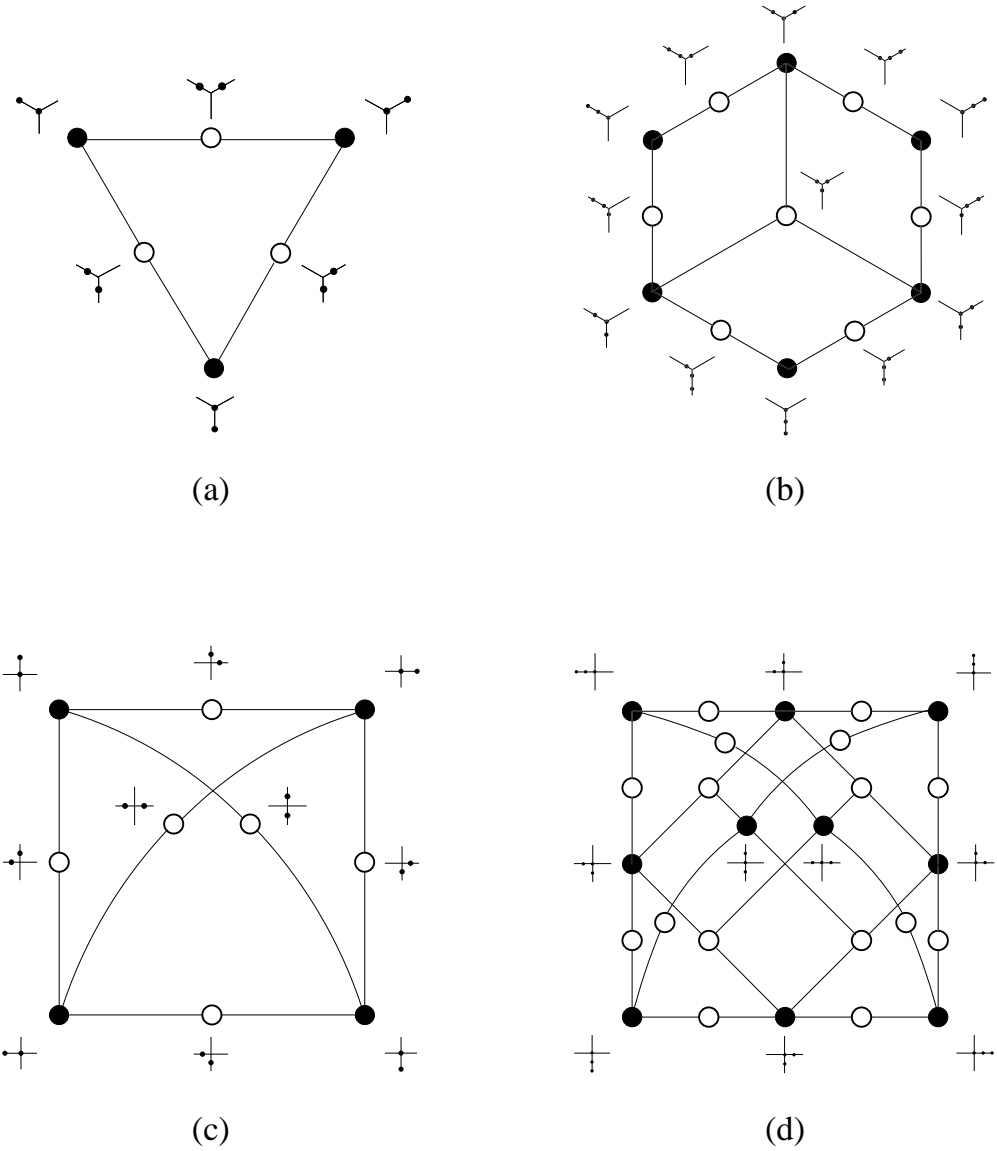
\begin{figure}
\begin{center}
\input{dksn.pstex_t}
\caption{Examples of the deformation retracts of (a) $D_2(S_3)$,
(b) $D_3(S_3)$, (c) $D_2(S_4)$, and (d) $D_3(S_4)$. The braid
groups are free on, respectively, one, three, three, and eleven
generators.  Solid vertices represent type I vertices and unfilled
represent type II. The small diagrams show the corresponding
configurations of points in $S_n$.} \label{dksn}
\end{center}
\end{figure}

\begin{ex}
By Corollary~\ref{group}, \[B_k(S_3) \cong
\mathcal{F}_{\binom{k}{2}},~B_2(S_n) \cong
\mathcal{F}_{\binom{n-1}{2}},~\mathrm{and}~B_3(S_n) \cong
\mathcal{F}_{\frac{2n+3}{3} \binom{n-1}{2}}\] where
$\mathcal{F}_p$ is the free group on $p$ letters.  See
Figure~\ref{dksn} for the illustrations of $D_k(S_n)$ as generated
by the method in the proof of Lemma~\ref{graph} in these four
cases.
\end{ex}

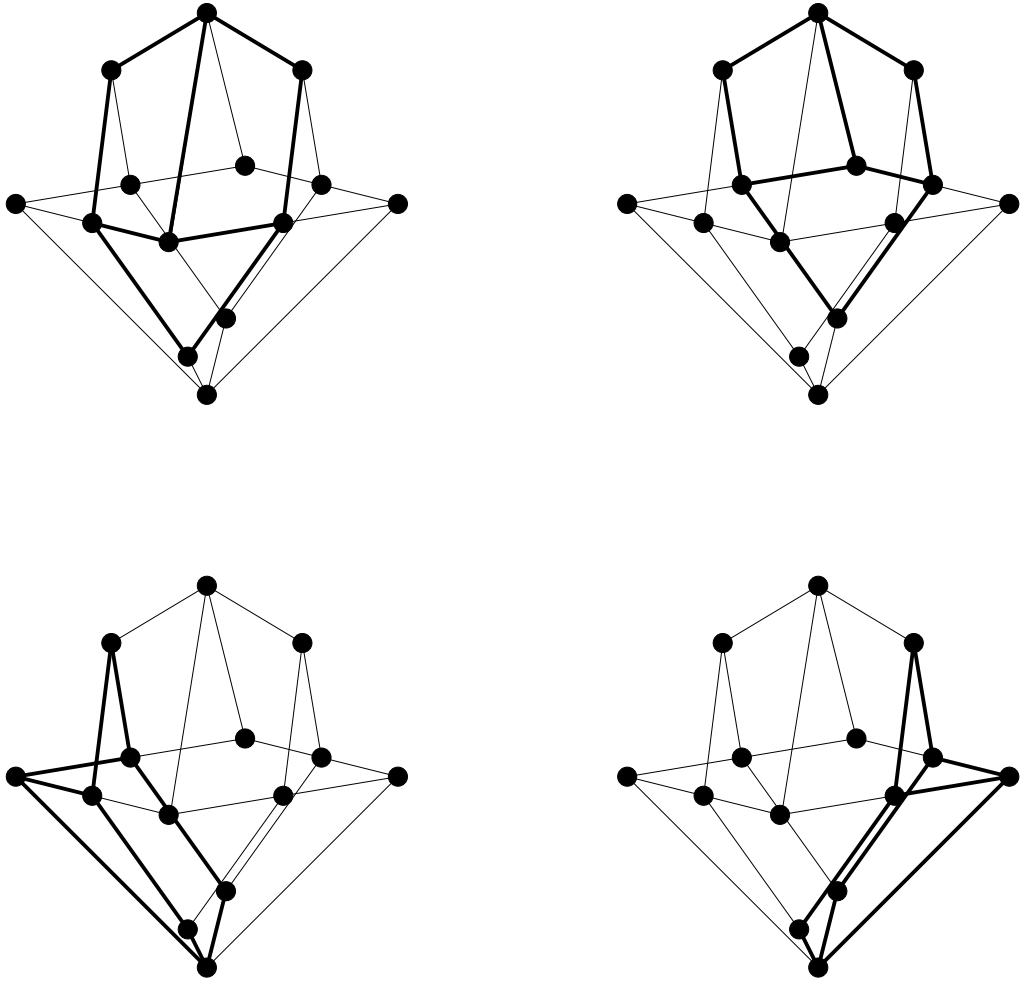
\begin{figure}
\begin{center}
\input{embeddings.pstex_t}
\caption{Four homeomorphic embeddings of $D_2(S_4)$, denoted by
solid edges, into $D_3(S_4)$.}\label{embeddings}
\end{center}
\end{figure}

Finally, Theorem~\ref{def} suggests a monomorphism which maps a
basis of $B_{k-1}(S_n,c_{k-1})$ into a basis of $B_k(S_n,c_k)$:

\begin{cor} \label{gen}
Let $\beta_j$ be a free basis for $B_j(S_n,c_j)$.  Then the
standard simplicial map $i_k:D_{k-1}(S_n) \hookrightarrow
D_K(S_n)$ induces a monomorphism $i_{k*}:B_{k-1}(S_n,c_{k-1})
\hookrightarrow B_k(S_n,c_k)$ which maps $\beta_{k-1}$ injectively
into $\beta_k$.
\end{cor}
\begin{proof}
Whenever there is a simplicial inclusion $i:L \hookrightarrow M$
for $L$ and $M$ one-dimensional complexes, their fundamental
groups $\pi_1(L)$ and $\pi_1(M)$ are free with bases $\beta(L)$
and $\beta(M)$ such that $i_*:\pi_1(L) \hookrightarrow \pi_1(M)$
is injective and $i_*(\beta(L)) \subseteq \beta(M)$.  Then
$D_{k-1}(S_n)$ and $D_k(S_n)$ are one-dimensional simplicial
complexes with a natural embedding $i_k:D_{k-1}(S_n)
\hookrightarrow D_k(S_n)$ which maps $c \in D_{k-1}(S_n)$ to $c
\cup \{x\}$ for $x$ the unique point such that $c \cup \{x\} \in
D_k(S_n)$.  (\emph{N.B.:} For a fixed $k$, such a point always
exists by our definition of the metric on $S_n$ and of
$D_j(S_n)$.) See Figure~\ref{embeddings} for four homeomorphic
embeddings of $D_{k-1}(S_n)$ into $D_k(S_n)$.
\end{proof}

\section{Conclusion}

In this paper, we demonstrate a one-dimensional deformation
retract of the configuration space of a star, $U_k^{top}(S_n)$.
This allows us to calculate the braid group $B_k(S_n,c_k)$ of the
star and gives a natural inclusion of braid groups
$B_{k-1}(S_n,c_{k-1}) \hookrightarrow B_k(S_n,c_k)$ which maps a
basis to a basis.  In a subsequent paper, the author and Frank
Connolly prove that, for any linear tree $T$, $B_k(T,c_k)$ is a
right-angled Artin group.

\section{Acknowledgements}

This paper was funded by the Provost's Office, Mathematics
Department, and Arts and Letters/Science Honors Program, all of
the University of Notre Dame.  Many thanks are owed to Frank
Connolly for his encouragement, support, and constant help.

\end{document}

%% file: vertices.pstex_t
\begin{picture}(0,0)%
\includegraphics{vertices.pstex}%
\end{picture}%
\setlength{\unitlength}{2763sp}%
\begingroup\makeatletter\ifx\SetFigFont\undefined%
\gdef\SetFigFont#1#2#3#4#5{%
  \reset@font\fontsize{#1}{#2pt}%
  \fontfamily{#3}\fontseries{#4}\fontshape{#5}%
  \selectfont}%
\fi\endgroup%
\begin{picture}(7667,4186)(-4124,-3362)
\end{picture}%

%% file: moves.pstex_t
\begin{picture}(0,0)%
\includegraphics{moves.pstex}%
\end{picture}%
\setlength{\unitlength}{2763sp}%
\begingroup\makeatletter\ifx\SetFigFont\undefined%
\gdef\SetFigFont#1#2#3#4#5{%
  \reset@font\fontsize{#1}{#2pt}%
  \fontfamily{#3}\fontseries{#4}\fontshape{#5}%
  \selectfont}%
\fi\endgroup%
\begin{picture}(7667,8986)(-4124,-8162)
\end{picture}%

%% file: dksn.pstex_t
\begin{picture}(0,0)%
\includegraphics{dksn.pstex}%
\end{picture}%
\setlength{\unitlength}{3158sp}%
\begingroup\makeatletter\ifx\SetFigFont\undefined%
\gdef\SetFigFont#1#2#3#4#5{%
  \reset@font\fontsize{#1}{#2pt}%
  \fontfamily{#3}\fontseries{#4}\fontshape{#5}%
  \selectfont}%
\fi\endgroup%
\begin{picture}(7808,8999)(-4195,-8912)
\end{picture}%

%% file: embeddings.pstex_t
\begin{picture}(0,0)%
\includegraphics{embeddings.pstex}%
\end{picture}%
\setlength{\unitlength}{3158sp}%
\begingroup\makeatletter\ifx\SetFigFont\undefined%
\gdef\SetFigFont#1#2#3#4#5{%
  \reset@font\fontsize{#1}{#2pt}%
  \fontfamily{#3}\fontseries{#4}\fontshape{#5}%
  \selectfont}%
\fi\endgroup%
\begin{picture}(7966,7666)(-382,-8544)
\end{picture}%